\documentclass[leqno]{article}

\usepackage{amssymb,amsmath,amsfonts,amsthm,mathrsfs}

\usepackage{epsfig}

\usepackage{pst-all}
\usepackage{multido}

\psset{xunit=1cm,yunit=1cm,runit=1cm}

\newcommand{\Con}{\ensuremath{\mathcal{C}}}
\newcommand{\Cinf}{\ensuremath{\mathcal{C}^\infty}}

\newcommand{\D}{\ensuremath{{\cal D}}}
\renewcommand{\S}{\mathscr{S}}

\newcommand{\Zin}{\ensuremath{C_*}}

\newcommand{\loc}{\ensuremath{\text{loc}}}

\newcommand{\mb}[1]{\ensuremath{\mathbb{#1}}}
\newcommand{\N}{\mb{N}}

\newcommand{\R}{\mb{R}}

\renewcommand{\d}{\ensuremath{\partial}}
\newcommand{\pdiff}[1]{\frac{\d}{\d#1}}

\newcommand{\grad}{\ensuremath{\mbox{\rm grad}\,}}
\renewcommand{\div}{\ensuremath{\mbox{\rm div}\,}}

\newfont{\bl}{msbm10 scaled \magstep2}

\newtheorem{theorem}{Theorem}[section]
\newtheorem{lemma}[theorem]{Lemma}
\newtheorem{proposition}[theorem]{Proposition}

\newtheorem{assumption}{Assumption}
\theoremstyle{definition}
\newtheorem{remark}[theorem]{Remark}
\newtheorem{example}[theorem]{Example}

\newcommand{\beq}{\begin{equation}}
\newcommand{\eeq}{\end{equation}}
\newcommand{\FT}[1]{\widehat{#1}}

\newcommand{\F}{\ensuremath{{\cal F}}}

\newcommand{\dis}[2]{\langle #1 , #2 \rangle}
\newcommand{\inp}[2]{\langle #1 | #2 \rangle}  %math mode
\newcommand{\notmid}{\mid\kern-0.5em\not\kern0.5em}

\newcommand{\norm}[2]{{\| #1 \|}_{#2}}

\newcommand{\ltwo}[1]{\norm{#1}{L^2}}

\newcommand{\al}{\alpha}
\newcommand{\be}{\beta}

\newcommand{\de}{\delta}
\newcommand{\eps}{\varepsilon}

\newcommand{\la}{\lambda}

\newcommand{\sig}{\sigma}

\newcommand{\inv}{^{-1}}

\newcommand{\SW}{\ensuremath{{\cal V}}}
\newcommand{\wt}[1]{\ensuremath{\widetilde{#1}}}

\begin{document}

\title{{\bf Evolution systems %and Sobolev-regularity
for paraxial wave equations of Schr\"odinger-type with non-smooth 
coefficients}\footnote{Work supported by FWF grants
P16820-N04 and Y237-N13}}

\author{\emph{Maarten de Hoop} \\
Center for Computational and
Applied Mathematics\\
Purdue University, West Lafayette, USA\\
\texttt{mdhoop@math.purdue.edu}\\
\ \\
\emph{G\"{u}nther H\"{o}rmann}\\
Fakult\"at f\"ur Mathematik\\
Universit\"at Wien, Austria \\
\texttt{guenther.hoermann@univie.ac.at}\\
\ \\
\emph{Michael Oberguggenberger}\\
Institut f\"ur Grundlagen der Bauingenieurwissenschaften\\
Universit\"at Innsbruck, Austria\\
\texttt{michael.oberguggenberger@uibk.ac.at}
}
\date{\today}
\maketitle

\vspace{-1cm}
\begin{abstract}%\enlargethispage{1cm}
  We prove existence of strongly continuous evolution systems in $L^2$for Schr\"odinger-type equations with non-Lipschitz coefficients in the principal part. The underlying operator structure is motivated
 from models of paraxial approximations of wave propagation in
  geophysics. Thus, the evolution direction is a spatial coordinate
  (depth) with additional pseudodifferential terms in time and low
  regularity in the lateral space variables. We formulate and analyze the
  Cauchy problem in distribution spaces with mixed regularity. The key
  point in the evolution system construction is an elliptic regularity
  result, which enables us to precisely determine the common domain of
  the generators. The construction of a solution with low regularity in
  the coefficients is the basis for an inverse analysis which allows to
  infer the lack of lateral regularity in the medium from measured
  data.
%\emph{AMS 2000 subject classification:}
\end{abstract}

%\setcounter{section}{-1}
%\input{intro}

%\section{Physical model}

\section{Introduction}

The paraxial equations in models of wave propagation are based on
parabolic symbol approximations in theories of wave operators. They
have been extensively applied in integrated optics, underwater
acoustic tomography as well as reflection seismic imaging (cf.\
\cite{Claerbout:85,WMcC:83}). They have also entered the analysis of
time-reversal mirror experiments with waves taking into account
stochastic variations in the wave speed (cf.\ \cite{BPR:02}). The
paraxial equations can be derived from the reduced wave or Helmholtz
equation and, since they split the wave fields according to a
prescribed principal direction of propagation, are also called one-way
wave equations. In particular, the leading-order parabolic symbol
approximation, also referred to as the narrow-angle or
beam-propagation approximation, leads to model equations of
Schr\"odinger-type. The well-posedness of the one-way wave Cauchy
problems has been discussed by \cite{HT:88,TH:86}. The methodologies
developed to date, however, have assumed smoothness of the wave speed
function (i.e., the coefficients in the wave operators).

In the analysis presented here, we depart from this smoothness
assumption by allowing the coefficients to be of any H\"older
regularity between zero and one, but typically non-Lipschitz. The
existence of distributional solutions to second order strictly
hyperbolic equations in general may fail below Log-Lipschitz
regularity of the coefficients (cf.\ \cite{CL:95}). In case of
H\"older regularity $2$ or higher a constructive approach for hyperbolic evolution equations has been developed in \cite{AHSU:07}, thereby extending results on propagation of singularities. The particular
equation considered here is derived from such a second order equation,
but the existence of its solution is not restricted by the same
conditions. Indeed, we exploit the framework of Sobolev space
techniques, in particular, multiplication of distributions in scales of
Sobolev spaces, to construct a strongly continuous evolution system.
The novelty of the paper lies in the method of construction which
not only provides a solution concept for the paraxial wave equation with
low coefficient regularity, but also allows us to investigate
how the coefficient regularity influences the solution.

\begin{figure}[ht]%
%\captionstyle{hang}%
%\psset{xunit=1cm,yunit=1cm,runit=1cm}
\begin{pspicture}(0,1)(12,-5)
\rput[tl](0,1){%
\centering\includegraphics[scale=0.7,angle=-90]{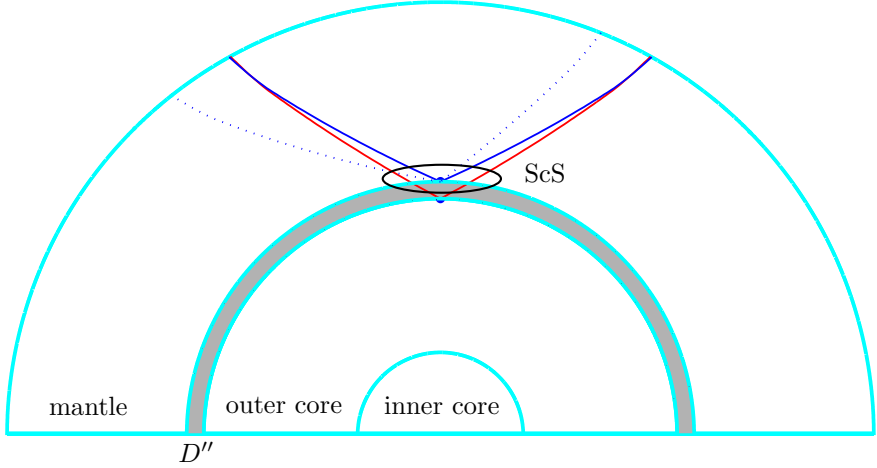}}
%\psgrid(0,1)(12,-5)
\psellipse(5.9,-1.5)(0.8,0.2)
\rput[tl](7,-1.3){ScS}
\rput[t](1.2,-4.4){mantle}
\rput[tl](2.4,-5){$D''$}
\rput[t](3.8,-4.4){outer core}
\rput[t](5.9,-4.4){inner core}
\end{pspicture}
\caption{The so-called $D''$ layer (in gray) above the core-mantle
  boundary (at approximately 2800km depth) with a core-reflected (ScS)
  wave and two precursors.  The precursors in the neighborhood of the
  top of $D''$ are locally modeled by our equation (indicated by
  ellipsoidal region), $x$ coinciding
  locally with the radial direction and $z$ coinciding with the
  tangential direction.}
\label{fig:ScS}%
\end{figure}

The class of coefficients of H\"older regularity between zero and one
arises in a variety of geophysical applications. Perhaps the most
fundamental one concerns the study of thermo-chemical boundary layers
and phase transitions in Earth's lowermost mantle --- the so-called $D''$ layer overlaying the core-mantle boundary (see Figure \ref{fig:ScS}; cf.\ \cite{science:07} for recent images of the phase transformation in $D''$).
Such phase transitions can only be probed by earthquake generated
seismic waves through scattering off these.  The relevant scattered
wave constituents appear as precursors to, for example, the
core-reflected compressional PcP phase and the horizontally polarized
shear ScS phase. The scattering is most prominent at large opening
angles (towards grazing incidence).
It is to this situation that the paraxial approximation and
coefficient dependence considered here applies.
In this context the principal direction, $z$, of propagation is perpendicular to
the direction of backscattering, $x$.
The central question is whether an imprint of the
coefficient regularity on the (regularity of the) scattered wave occurs.
A positive answer to this question, as provided in this paper, allows to infer a lack of lateral regularity in the medium from a lack of regularity of the measured data.

Seismic inverse scattering has been formulated mathematically in terms of evolution equations with respect to the depth variable in cf.\ \cite{SdH:02,SdH:05,SdH:06} (with smooth symbols in the single scattering approximation). At the basis of these models are one-way wave equations, which are typically of the form
$$
  \big( \d_z \pm i B(z,x,D_t,D_x) \big) u = f,
  $$
 where $B$ is (microlocally) a pseudodifferential operator with principal part
$$
   b(z,x,\tau,\xi) =\tau c(z,x)^{-1} \sqrt{1 - \tau^{-2} c(z,x)^2
   |\xi|^2}.
$$
Here, $t\in\R$ is time, $z\geq 0$ denotes depth, $x \in \R^d$
are lateral spatial directions, and $D = i\inv \d$. Approximation of the
square root to leading order results in the standard
Schr\"odinger-type paraxial equation
\[
   \big( \d_z + \frac{i}{c(z,x)} D_t \big) w + \frac{1}{2 i} D_t\inv
   \sum\limits_{j=1}^{d} D_{x_j} \big( c(z,x) D_{x_j} w \big) = 0,
\]
where $c(z,x)$ is the local speed of propagation and $z$ plays the
role of the evolution parameter. In the frequency
($\tau$-)domain the above equation is transformed with a so-called co-moving
frame of reference according to $\FT{w}(z,x,\tau) =
\FT{u}(z,x,\tau)$ $\exp({i \tau T(z,x)})$, where
$ T(z,x) = \int_0^z dz' / c(z',x)$.
Then the paraxial equation attains the form
$$
  \big( \d_z - i \bar{A} \big) \FT{u} = 0
$$
in which $\bar{A}$ is given by
\[
   \bar{A} \FT{u} = \frac{1}{2} \sum\limits_{j=1}^{d}
   e^{i \tau T(z,x)}
    D_{x_j}\Big( c(z,x) \tau\inv D_{x_j}\big(
    e^{- i \tau T(z,x)} \FT{u} \big)
    \Big).
\]
The second-order differential operator $\bar{A}$
can be written as the sum of the self-adjoint operator
\[
   A(\tau;z,x,D_x) := \frac{1}{2}
   \sum\limits_{j=1}^{d} D_{x_j} c(z,x) \tau\inv  D_{x_j}
\]
and a symmetric perturbation. We observe that if $\norm{D_x
  T(z,.)}{L^\infty(\R^d)} < K$, where $K$ is an appropriate constant,
then it is guaranteed that this perturbation is $A$-bounded with relative
bound less than $1$ (cf.\ \cite[Section X.2]{RS:V2}). From the
viewpoint of generators of strongly continuous contraction semigroups on
$L^2(\R^d)$ the simplification of $\bar{A}$ by $A$ is of no consequence.

\begin{remark}[{\bf Directional decomposition and one-way wave equations}]

Directional decomposition leads to the introduction of one-way wave
equations \cite{taylor1975}. One-way wave equations approximate
solutions to the wave equation microlocally, relative to a principal
direction of propagation. (This principal direction does not need to
be defined globally; one can introduce curvilinear coordinates and an
associated Riemannian metric to generate such directions locally.) The
validity of one-way wave propagation breaks down when singularities
tend to propagate in a direction perpendicular to the principal
direction (that is, in a transverse direction). Hence, to make the
statement concerning approximation above, precise, one needs to
introduce a microlocal attenuation \cite{stolk2004}. The mentioned
procedures and results require smooth coefficients and symbols, and
can be proven by making use of the calculus of pseudodifferential
operators and Fourier integral operators with complex phase.

It has been demonstrated that, in special cases, one can weaken the
condition of smooth coefficients. For example, if the coefficients are
independent of the coordinate along the principal direction, one can
allow a step function (in a transverse coordinate) and still solve the
associated scattering problem by methods of one-way wave
equations. The approach to carry out such an evaluation can be found
in \cite{fishman2000}. Indeed, scattering in the transverse directions
can, at least in special cases, be incorporated in the one-way wave
equation. Moreover, in the case of wave propagation in random media,
the (stochastic) paraxial equation naturally appears as well
\cite{papanicolaou2004}. However, a general result concerning
directional decomposition for -- or recomposition to solutions of --
the wave equation with non-smooth coefficients has not been obtained.

In this paper, we address the general problem of ``transverse
scattering'' by a one-way wave equation, which constitutes one
component in the development of a general theory referred to above.
\end{remark}

\begin{remark}[{\bf Regularity and the second order wave equation}]
Both, mode decoupling of a second order wave equation into one-way wave equations
as well as the derivation of the narrow beam approximation outlined above, require higher order differentiablity of the coefficient $c(z,x)$ with respect to all variables to make sense (due to truncation of symbol expansions). However, the resulting paraxial wave equation displays precisely the same coefficient regularity as in the original second order wave equation. Hence as a model equation it still reflects the correct medium properties on all scales. The (exact) solutions of the paraxial equation then serve as a narrow beam approximation to solutions of the original wave equation. In particular, the regularity properties are comparable on the same scales.

The fine tuned well-posedness theorem for wave equations by Colombini and Lerner in \cite{CL:95} assumes Log-Lipschitz regularity of the coefficients in the principal part. Moreover, their results are sharp in the sense that counterexamples to solvability exist when the coefficient regularity is below Log-Lipschitz (but still of any continuity type arbitrarily close to such). 
In order to indicate how their results relate to ours, we repeat the key energy estimate from \cite{CL:95}:

A function $a \in L^\infty{(\R^d)}$ is said to be a \emph{Log-Lipschitz} function if
$$
  \norm{a}{LL} :=
   \sup_{x \in \R^d} |a(x)| +
   \sup_{x \neq y \in \R^d, \atop |x-y| \leq 1/2}
   \frac{|a(x) - a(y)|}{|(x - y) \log(|x-y|)|} < \infty.
$$
Colombini and Lerner consider second order wave operators of the form
$$
   Pu := \d_t^2 u - \sum_{1 \leq j,k\leq n} \d_{x_j} \big(a_{jk}(x,t) \d_{x_k} u \big) + M(x,t,\d_t,\d_x) u,
$$
where $(a_{jk})_{1 \leq j,k\leq n}$ is real symmetric with Log-Lipschitz components and satisfies with some $1 \geq \de_0 > 0$ the strong ellipticity condition
$$
   \sum_{1 \leq j,k\leq n} a_{jk}(x,t) \xi_j \xi_k
   \geq \de_0 |\xi|^2
   \qquad  (\xi \in \R^n, (x,t)\in \R^{n+1}),
$$
and $M(x,t,\d_t,\d_x)$ is a first-order differential operator with H\"older-continuous coefficients.

Let $P$ be as above and $\theta \in ]0,1/4]$. There exist $\be > 0$ and $C > 0$ such that for $u \in \Cinf(\R^{n+1})$ and $t \in [0,1/\be]$ the \emph{energy estimate} 
\begin{multline*}
  \sup_{0 \leq s \leq t}
  \norm{\d_t u(.,s)}{H^{-\theta - \be s}} +
  \sup_{0 \leq s \leq t}
  \norm{u(.,s)}{H^{1-\theta - \be s}}\\
  \leq C \left(
   \int_0^t \norm{Pu(.,s)}{H^{-\theta - \be s}} \, ds
   + \norm{\d_t u(.,0)}{H^{-\theta}} +
   \norm{u(.,0)}{H^{1-\theta}}
  \right),
\end{multline*}
holds (cf.\ \cite[Equation (2.6)]{CL:95}), where $\be$ depends only on $\de_0$, on the Log-Lipschitz norm of the $a_{jk}$, and on the H\"older norms of the coefficients in $M(x,t,\d_t,\d_x)$.

We show in the sequel that for the paraxial wave equation the condition on the coefficient regularity can be relaxed. For example, if $\eps > 0$ any function in $H^{1+\eps}(\R^2)$ of local behavior like $x \mapsto |x|^{1/2 + \eps}$ is not Log-Lipschitz continuous but satsifies the assumptions of our main results below.

The low coefficient regularity in our model conditions has its price in terms of a few technical aspects of the current paper: Additional care is needed in identifying the appropriate distribution and function spaces that allow for the description of mixed regularity properties and for the rigorous formulation of a solution concept. The existence proof then consists in showing a series of functional analytic properties to establish an evolution system of operators; among these the basic self-adjointness property --- in the disguise of an elliptic regularity lemma --- is derived by employing rather delicate regularity properties of multiplication in certain subspaces of the space of distributions.
\end{remark}

The plan of the paper is as follows. In Section 2 we present the
precise form of the operator and specify our (low) regularity
assumptions on the coefficients. The solution will be sought as a
continuous map of depth into the space of temperate $L^2$-valued
distributions. Section 3 is devoted to the construction of the
evolution system in the frequency domain. First, we prove that
$A(\tau;z,x,D_x)$ generates a unitary group at fixed $z$ and $\tau$.
The determination of its domain requires delicate use of the duality
product of distributions as well as a bootstrap argument involving
multiplication in scales of Sobolev spaces. This leads to the
construction of an evolution system at fixed $\tau$. Finally, strongly
continuous dependence on the frequency parameter $\tau$ is
established, based on a difference approximation. Again a subtle
interplay between regularity arguments and distributional products is
at the heart of our arguments. The strong continuity enables us, in
Section 4, to construct a solution of the evolution system in
frequency domain with distributional data. Finally, existence,
uniqueness, and regularity of solutions to the original Cauchy problem
is obtained. As an application to inverse regularity analysis, we
obtain that a lack of $H^2$-regularity in the observed solution
implies the existence of a region in which the lateral regularity of
the medium is at most $H^1$ on the Sobolev scale.

%%% Local Variables:
%%% mode: latex
%%% TeX-master: "paraxial"
%%% End:

\section{The Cauchy problem: function spaces and \\ coefficient regularity}

We recall the definition of temperate distributions on $\R$ with values
in a Banach space $E$ (cf.\ \cite[Section 39.3]{Treves:75};
but note that we use a different topology here):
let $\S'(\R;E)$ be the space of
continuous linear maps $\S(\R) \to E$, equipped with the topology of pointwise
convergence; the Fourier transform $\F$ on $\S(\R)$ is extended to $\S'(\R;E)$
by setting $(\F G)(\phi) = G(\F \phi)$ ($G \in \S'(\R;E)$, $\phi\in\S(\R)$),
which is easily seen to be an isomorphism (for the locally convex structure).

We denote the time variable by $t\in\R$ and introduce coordinates $z
\in [0,\infty)$ for depth (the evolution direction in our context) and
$x\in\R^d$ for the lateral variation, where $d \leq 2$. As basic space
of the \emph{wave components} we consider \beq\label{basic_space} \SW
:= \Con([0,\infty), \S'(\R;L^2(\R^d)).  \eeq Its elements are
continuous maps $z \mapsto u(z)$ of the depth variable $z$ into
temperate distributions of time $t$ valued in $L^2$-functions of the
lateral variables $x$.  When we need to keep track of precise
regularity information in the lateral variation of the waves, we may
employ the Sobolev-scale $H^s(\R^d)$ ($s\in \R$) and define
\beq\label{basic_s_spaces} \SW^s := \Con([0,\infty), \S'(\R;
H^s(\R^d)).  \eeq

Let $\F_t : \S'(\R; H^s(\R^d)) \to \S'(\R; H^s(\R^d))$ denote (partial) Fourier
transform with respect to the time variable. We extend $\F_t$ to an
isomorphism $\wt{\F_t}$ of (the locally convex structure of) $\SW^s$ by
\[
    (\wt{\F_t} u)(z) := \F_t (u(z,.)) \qquad \forall z \geq 0.
\]

We consider the following Cauchy problem for a prospective solution
$u \in \SW$ ($\div$ and $\grad$ with respect to $x\in\R^d$):
\begin{align}
    P u := \d_z u - i\; \div( C(z,x,D_t)\cdot \grad u) &= f \in \SW
                                                \label{orig_PDE}\\
                               u \mid_{z = 0} &= u_0 \in \S'(\R;L^2(\R^d)).
                                                \label{orig_initial_cond}
\end{align}
Here, $C$ is a pseudodifferential operator in $t$  with
parameters $z$ and $x$.  While in the classical paraxial wave equation it is
of order $1$, here we may assume it is of some order  $m\in\R$.
The precise conditions are collected in the  following
\begin{assumption}
The symbol of $C(z,x,D_t)$ is of the form
\beq\label{C_structure}
    C(z,x,\tau) = c(z,x,\tau) \cdot I_d =
        \Big( c_0 + \sum_{l=1}^N c_l(z,x) h_l(\tau) \Big) \cdot I_d,
\eeq
where $N\in\N$, $I_d$ is the $d\times d$ identity matrix and the following
hold:
\begin{enumerate}
 \item For $l = 1,\ldots,N$: $h_l$ is a real-valued smooth symbol (of order $m$) on
 $\R$, i.e., for all $k\in\N_0$ an estimate $|\d_\tau^k h_l(\tau)| = O(|\tau|^{m-k})$
 holds when $|\tau|$ is large; in addition, we assume that
 $|h_l(\tau)| \geq \eta_0$  near $\tau = 0$ ($l = 1,\ldots,N$) with
 some constant $\eta_0 > 0$ (this can be achieved by adding a cut-off
 function without changing the relevant frequency range).
 \item $c_0$ is a positive constant.
 \item there is an $r\in (0,1)$ such that for $1 \leq l \leq N$:
      $c_l$ is in $\Con^1([0,\infty),H^{r + 1}(\R^d))$ and real-valued.
 \item for all $(z,x,\tau) \in [0,\infty)\times\R^d\times\R$:\enspace
    $c(z,x,\tau) \geq c_0$.
\end{enumerate}
\end{assumption}
\begin{remark}
(i) The operator action corresponding to a typical term $c_l(z,x) h_l(D_t)$ in
the sum decomposition (\ref{C_structure}) on any element $w$ of $\SW$ is given
as follows: for all $\phi\in\S(\R)$
\[
    \big(h_l(D_t) w(z)\big)(\phi) = w(z) (\F(h_l\,\F\inv\phi)) \in L^2(\R^d),
\]
which is then multiplied by the function $c_l(z,.)$.

(ii) Note that parts (i-iii) of the Assumption imply, for any $z$ and $\tau$
fixed, the following (Zygmund-)H\"older-continuity:
\[
    c(z,.,\tau) - c_0 \in H^{r+1}(\R^d)
            \subseteq \Zin^{r + 1 - \frac{d}{2}}(\R^d)
\]
(cf.\ \cite[Proposition 8.6.10]{Hoermander:97}). Thus the coefficients
have lateral H\"older-regularity of order $r + 1 - \frac{d}{2}$. In
the most relevant case from geophysics, $d = 2$, this yields
coefficients in $\Zin^r(\R^2)$, but not necessarily Lipschitz
continuous. (This includes the situation of a boundary layer, as
discussed in the introduction, where the coefficient is smooth in one
of the two variables.) Part (iv) implies uniform ellipticity of the
lateral differential operator.

(iii) Note that in dimension $d \geq 3$ Sobolev regularity of order $r
+ 1$ would not imply continuity of the coefficients (in case $1 < r + 1 <
d/2$).
\end{remark}

\begin{example}
We consider a model with two-dimensional lateral variation ($d=2$) in
the medium properties of low H\"older regularity depending on depth.
In (\ref{C_structure}) we put $c_l = 0$ when $l \geq 2$ and let $c_1$ be
of the form
\[ 
 c_1(z,x) = \chi(z,x)|x|^{\al(z)}
\] 
with the following specifications:
%\begin{trivlist} 
$\chi \in \Con^1([0,\infty)\times \R^d)$ such that $\chi(z,x) =
\chi_0(z)$ 
when $|x| \leq R_1$ and $\chi(z,x) = 0$ when $|x| \geq R_2$ for
certain radii $0 < R_1 < R_2$ and some $\chi_0 \in \Con^1([0,\infty))$;
$\al \in \Con^1([0,\infty))$ with some uniform positive lower
  bound $\al_0$, i.e.,  $\al(z) \geq \al_0 > 0$ for all $z$.
%\end{trivlist}
Then we may choose any $r$ such that $0 < r < \al_0$ and obtain the following
regularity properties at arbitrary fixed values of $z$ and $\tau$:
\[
    c(z,.,\tau) - c_0 \in H^{r+1}(\R^2)
            \cap \Zin^{\al(z)}(\R^2) \subseteq  H^{r+1}(\R^2)
            \cap \Zin^{\al_0}(\R^2) .
\]
(Observe that locally in two dimensions, for any $0 < s < 1$,
the function $|x|^s$ belongs to $\Zin^s$ and to $H^{s + 1 -
  \eps}$ for every $\eps > 0$ but not to $H^{s+1}$.)
\end{example}

Applying $\wt{\F_t}$ to (\ref{orig_PDE}-\ref{orig_initial_cond}) we
obtain an equivalent formulation of the Cauchy problem 
in the frequency domain:
\begin{align}
    \wt{P} v = \d_z v - i\; \div( c(z,x,\tau)\, \grad v) &= g \in \SW
                                                        \label{PDE}\\
                               v \mid_{z = 0} &= v_0 \in \S'(\R;L^2(\R^d)).
                                                \label{initial_cond}
\end{align}
Equation (\ref{PDE}) is an evolution equation for depth $z$ with the
second-order operator
 \beq \label{A_operator}
    A(\tau;z,x,D_x) v := \div( c(z,x,\tau)\, \grad v)
 \eeq
 acting in the lateral $x$-domain and smoothly
depending on the ``external'' parameter $\tau$. Note that
$A(\tau;z,x,D_x)$ is uniformly elliptic by Assumption 1 (iv).

\begin{remark}\label{damping_remark} Note that $P$ in (\ref{orig_PDE})
commutes with convolution in the time variable.
Therefore, damping
(or cut-off) of high frequencies in the data of (\ref{PDE}-\ref{initial_cond})
corresponds to time-smoothing the data in the original problem (\ref{orig_PDE}-\ref{orig_initial_cond}): more precisely, if
a frequency filter $\FT{\chi} \in \S(\R)$ is applied by $g := \wt{\F_t}(f)
\cdot \FT{\chi}(\tau)$, $v_0 := \F_t (u_0) \cdot \FT{\chi}(\tau)$ and
$v$ is a solution to (\ref{PDE}-\ref{initial_cond}) then $u := \wt{\F_t}\inv
(v)$ solves
(\ref{orig_PDE}-\ref{orig_initial_cond})
with the data changed to $f \mathop{*}\limits_{(t)} \chi$ and
$u_0 \mathop{*}\limits_{(t)} \chi$.
\end{remark}

%%% Local Variables: 
%%% mode: latex
%%% TeX-master: "paraxial"
%%% End: 

\section{Evolution system}

In this section, we will show that, up to any finite depth $Z > 0$,
the family of unbounded operators $i \, A(\tau;z,x,D_x)$ ($\tau\in\R$,
$z \geq 0$) generates a strongly continuous evolution system (or
fundamental solution) $\{ U(\tau;z_1,z_2) : 0 \leq z_1 \leq z_2 \leq
Z\}$ on $L^2(\R^2)$ in the sense of \cite[Chapter 4]{Tanabe:79} which,
in addition, is strongly continuous in the frequency variable 
$\tau\in\R$. In a
first step we freeze both parameters, $\tau$ as well as $z$, and
construct a strongly continuous (semi-)group of operators on
$L^2(\R^d)$.

\subsection{Unitary group at frozen values of $\tau$ and $z$}

\paragraph{Notational simplifications:} By abuse of
notation we will employ the short-hand symbols $A := A(\tau;z,x,D_x)$, $c(x) :=
c(z,x,\tau)$, and $c_1(x)$ now denoting $\sum_{l \geq 1} c_l(z,x) h_l(\tau)$.
%Furthermore, in Assumption 1 we may restrict to the case $\tau > 0$: in fact,
%if $c$ is negative the line of reasoning following below can be repeated with
%$-A$ in place of $A$; on the other hand, $c \equiv 0$ yields $i \, A = 0$,
%which is the generator of the trivial group.
To summarize, using the above conventions and Assumption 1, we have
 \beq\label{A_defined}
    A v = \div( c(x)\, \grad v)
 \eeq
 as unbounded, formally self-adjoint operator on $L^2(\R^d)$ with coefficient
 \beq\label{coeff_structure}
    c(x) = c_0 + c_1(x) \quad \text{ with } c_0 > 0
        \text{ and } 0 \leq c_1 \in  H^{r+1}(\R^d).
 \eeq
We will show that $A$ is a self-adjoint operator with domain $D(A) =
H^2(\R^d)$. 

\begin{remark}
Note that self-adjointness of $A$ with domain $H^2$ would be immediate
from uniform ellipticity in case the coefficient were smooth. On the other
hand, self-adjointness on some domain could be obtained in an abstract fashion via
quadratic forms under mere $L^\infty$-assumptions (\cite[Section
VIII.6] {RS:V1}).  However, in accordance with our
focus on the interplay of the coefficient regularity class with  
qualitative solution properties, we will give an explicit domain description in
terms of Sobolev spaces, which in addition is uniform with
respect to $\tau$ and $z$.  
\end{remark}

Observe that, due to the low coefficient regularity, we also have to establish that
$A$ is well-defined on all of $H^2(\R^d)$ as an operator into $L^2(\R^d)$. This
is included in the following lemma as the special case $s = 0$.
\begin{lemma}\label{inclusion_lemma} Let $0 \leq s < r < 1$ and $v \in
H^{s+2}(\R^d)$. Then $A v \in H^s(\R^d)$.
\end{lemma}
\begin{proof} Each component of $\grad v$ is in $H^{s+1}(\R^d)$. Hence
multiplying with the $H^{r+1}$-coefficient $c_1$, as well as with the constant
$c_0$, is well-defined within $H^{s+1}(\R^d)$ since this is an algebra.
\end{proof}
We set $D(A) := H^2(\R^d)$ and note that an integration by parts immediately
yields that $A$ is symmetric, i.e., $D(A) \subseteq D(A^*)$ and $A^*
\mid_{D(A)} = A$. We proceed to show that also $D(A^*) \subseteq D(A)$ by which
self-adjointness will be established.

By definition, the adjoint operator has domain
 \[
    D(A^*) = \{ v \in L^2 \mid \text{for some } w\in L^2: 
  \inp{\psi}{w} = \inp{A \psi}{v} \text{ for all } \psi \in H^2  \},
 \]
where $\inp{\ }{\ }$ denotes the inner product in $L^2(\R^d)$. Let $v \in
D(A^*)$ and choose a sequence $(v_k)_{k\in\N}$ in $H^2(\R^d)$ which converges
to $v$ in $L^2(\R^d)$. On the one hand, there exists $w\in L^2(\R^d)$ such that
we have for all test functions $\psi\in \D(\R^d)$
 \[
    \inp{A^* v_k}{\psi} = \inp{v_k}{A \psi} \to \inp{v}{A \psi} =
    \inp{w}{\psi} \quad \text{ when } k\to\infty.
 \]
Thus, $(A^* v_k)_{k\in\N}$ has the distributional limit $w \in L^2(\R^d)$. On
the other hand, since $c_1 \in H^1(\R^d)$ and $\grad v_k \to \grad v$ in
$H^{-1}(\R^d)$ (as $k\to\infty$) we may employ the continuous duality
product of distributions 
(cf.\ \cite[Proposition 5.2]{O:92}) $H^1 \times H^{-1} \to W^{-1,1}_\loc$ and
deduce that $A^* v_k = A v_k \to \div( c(x)\, \grad v) = Av$ in
$W^{-2,1}_\loc$, hence in the sense of distributions. By uniqueness of
distributional limits, we deduce that $Av = w \in L^2(\R^d)$. We obtain
 \[
    D(A^*) = \{ v \in L^2 \mid Av \in L^2 \}.
 \]
The assertion $D(A^*) \subseteq D(A) = H^2(\R^d)$ follows now from the
following result. For later reference, it is stated in slightly more general
terms. (We first consider the important case $d = 2$ and leave the case of
one-dimensional lateral variation for a remark below.)

\begin{lemma}[Elliptic regularity] \label{elliptic_regularity_prop} Let 
$0 \leq s < r < 1$ and $v \in H^{s}(\R^2)$ such that $A v \in H^s(\R^2)$. 
Then $v$ belongs to $H^{s+2}(\R^2)$.
\end{lemma}

The proof will be based on repeated use of the following three facts, which we
collect in a preparatory list of ``ingredients'':
\begin{description}
 \item[Fact A:] Let $w_j \in H^{s_j}(\R^2)$ ($j=1,2$) such that $s_1 + s_2
 \geq 0$. Then
 \[
    w_1 \cdot w_2 \in H^{s_0}(\R^2),
 \]
 where
 \[
    s_0 = \begin{cases}
            \min(s_1, s_2, s_1 + s_2 - 1) & \text{ if }
                            s_1 \not= \pm 1, s_2 \not= \pm 1,
                            \text{ and } s_1 + s_2 \not= 0\\
            \min(s_1, s_2, s_1 + s_2 - 1 - \eps) & \text{ with }
                \eps > 0 \text{ arbitrary, otherwise}.
    \end{cases}
 \]
 This is included in the statement of \cite[Theorem 8.3.1]{Hoermander:97}.
As can be seen from the proof therein, one also obtains
%the following estimate
% (with a constant $C_{s_1,s_2}$ independent of $w_1$ and $w_2$)
% \[
%   \norm{w_1 \cdot w_2}{H^{s_0}} \leq C_{s_1,s_2} \norm{w_1}{H^{s_1}}
%                                       \norm{w_2}{H^{s_2}},
% \]
%   which implies
 continuity of the multiplication $H^{s_1} \times H^{s_2} \to H^{s_0}$ (with
 respect to the corresponding Sobolev-norms).
 \item[Fact B:] We can find a function $F\in \Cinf(\R)$, $F(0) = 0$, such that
 \[
    \frac{1}{c(x)} = \frac{1}{c_0} + F(c_1(x)).
 \]
 In particular, we obtain $1/c - 1/c_0 = F(c_1) \in H^{r+1}(\R^2)$.
 Since $r > 0$, this follows from \cite[Theorem 8.5.1]{Hoermander:97}
 once $F$ is given. To find $F$, we simply set $F(y) := - y / (c_0 (c_0 + y))$
 when $y \geq c_0/2$ and extend it in a smooth way to $\R$ such that $F(0) = 0$.
 \item[Fact C:] If $v\in L^2(\R^2)$ the equation
  \[
    \tag{C1} Av = \div(c \, \grad v) = \Delta(c\, v) - \div(v \, \grad c)
  \]
   holds in $\D'(\R^2)$, where the occurring products are defined as follows:
   using $\grad v \in H^{-1}$ we get $ c \, \grad v \in W^{-1,1}_\loc$ by the
   duality product \cite[Proposition 5.2]{O:92}; $c \, v \in L^2$ since
   $c\in L^\infty$; and $v \, \grad c \in L^1$ because $\d_j c \in H^r \subseteq L^2$.\\
   Under the stronger assumption $v\in H^{r+1}(\R^2)$ we have, in addition,
   that
 \[
    \tag{C2} Av = \grad c \cdot \grad v +  c \, \Delta v
 \]
 in $\D'(\R^2)$, with the meaning of the products on the right-hand side as
 follows: since $\Delta v \in H^{r-1}$, Fact A applies and yields $c \Delta v \in
 H^{r-1}$; furthermore, $\grad c$ and $\grad v$ both lie in $H^r$, so
 that another application of Fact A shows that their (Euclidean inner) product
 belongs to $H^{\min(r,2r-1)}$.
\end{description}

\begin{remark} Note that formula (C2) represents $A$ as an operator with a
 (H\"older-) continuous coefficient in its principal part and Sobolev
 regularity in the lower orders. We observe that in such situation,
\cite[Theorem 17.1.1]{Hoermander:V3} gives local solvability in $H^2$ for
 right-hand sides in $L^2$. However, the latter does not imply
 $H^2$-regularity of any $L^2$-solution. For the pure regularity
 question, it also seems that methods based on perturbations of constant
 coefficient operators do not apply either, since $A$ is not
 necessarily of constant strength (cf.\ \cite[Chapter XIII]{Hoermander:V2}). 
\end{remark}

\begin{proof}[Proof of Lemma \ref{elliptic_regularity_prop}] To begin with,
we only know that $v$ as well as $Av$ belong to $H^s(\R^2)$. The proof proceeds
in three steps, successively revealing higher regularity.

\emph{Claim 1:} $v \in H^{s+r}(\R^2)$

We have $\grad c \in H^r$, so that application of Fact A, noting that $r-1 <
0$, gives $v \, \grad c \in H^{s+r-1}$, hence $\div(v \, \grad c) \in
H^{s+r-2}$. Since $Av \in H^s$ we deduce from equation (C1) that $\Delta(c v)
\in H^{s+r-2}$, which in turn implies that $c v \in H^{s+r}$. Now invoke the
decomposition from Fact B and write
\[
    v = \frac{1}{c_0} c v + F(c_1) c v.
\]
The first part clearly is in $H^{s+r}$, and for the second summand the same is
true by Fact A. (Note that $\min(s+r, s + 2r - \eps) = s+r$ if $0 < \eps < r$.)

\emph{Claim 2:} $v \in H^{r+1}(\R^2)$

We may start from $v \in H^{r + s}$ by claim 1 and proceed inductively to show
that
\[
   \tag{$\star$} v \in H^{r + \min(1, s + j r / 2)} \quad j \geq 0.
\]
Claim 2 then follows upon choosing $j$ sufficiently large (i.e., $j \geq
2(1-s)/r$ steps will be required).

The case $j=0$ is just claim 1, so we assume that the assertion holds for some
$j \geq 0$. If $s + j r / 2 \geq 1$ it is trivially satisfied for larger values
of $j$, therefore we assume $t_j := s + j r / 2 < 1$ and that
\[
    v \in H^{r + t_j}.
\]
Fact A gives $v \, \grad c \in H^{r+t_j} \cdot H^r \subseteq H^{\min(r,t_j + 2r
- 1 -\eps)} \subseteq H^{\min(r,t_j + 3 r / 2 - 1)}$ upon choosing $\eps <
r/2$. Thus, using the short-hand notation $r_j := \min(r+1,t_j + 3r/2)  \leq r
+ 1$ we may infer that $\div(v \, \grad c) \in H^{r_j - 2}$. Since $r_j - 2
\leq r - 1 < 0 < s$ we also have $Av \in H^s \subset H^{r_j - 2}$. Equation
(C1) now implies that $\Delta(c v) \in H^{r_j - 2}$, hence $c v \in H^{r_j}$.
Again by Fact B, combined with Fact A, we obtain
 \[
    v = \frac{1}{c_0} c v + F(c_1) c v \in H^{r+1} \cdot H^{r_j}
        \subseteq H^{r_j} = H^{r + \min(1,t_j + r/2)},
 \]
 which means ($\star$) for $j+1$ in place of $j$.

\emph{Claim 3:} $v \in H^{s+2}(\R^2)$

We use a similar strategy as in the proof of claim 2 and prove inductively that
\[
    \tag{$\star\star$} v \in H^{\min(s+2, 1 + (j+1)r / 2)} \quad j \geq 1.
\]
Claim 3 then follows when $j$ is chosen sufficiently large (i.e., $j \geq
2(1+s)/r -1$ steps are required).

The basic case $j=1$ corresponds to claim 2, so we proceed with some $j \geq
1$, under the additional assumption $r \leq s_j := (j+1) r / 2 < s + 1$ to
exclude trivial cases. In other words, we assume that
\[
    v \in H^{s_j + 1}.
\]
Therefore, again by Fact A and choosing a possibly occurring $\eps < r/2$, we
deduce
\[
    \grad c \cdot \grad v \in H^r \cdot H^{s_j} \subseteq
        H^{\min(r,s_j + r/2 - 1)} = H^{\min(r, q_j)},
\]
where we have introduced $q_j := s_j + r/2 - 1$. Exploiting equation (C2) and
noting that $Av \in H^s$, $s < r$, we extract the information that $c\, \Delta
v \in H^{\min(s,q_j)}$. Once again we use the decomposition corresponding to
Fact B and the statement of Fact A to deduce
\[
    \Delta v \in H^{\min(s,q_j)} + H^{r+1} \cdot H^{\min(s,q_j)} \subseteq
        H^{\min(s,q_j)}
\]
and a fortiori that $\Delta v \in H^{\min(s+2,q_j + 2)}$. But $q_j = (j+1) r /2
+ r/2 -1$ hence $q_j + 2 = 1 + (j+2) r / 2$ and ($\star\star$) is proven with
$j+1$ in place of $j$.
\end{proof}

\begin{remark} The one-dimensional analogue of Lemma
  \ref{elliptic_regularity_prop} is more elementary: $Av = (c v')' \in
  H^s(\R)$ implies $c v' \in H^{s+1}(\R)$ and, since Fact B is valid
  for $d=1$ as well, we obtain $v' \in H^{s+1}(\R)$; thus, $v \in
  H^{s+2}(\R)$.
\end{remark}

We summarize the intermediate conclusions from the discussion so far in a
separate statement, where we appeal to Stone's theorem providing us
with the exponential unitary group $T(z) = \exp(i z A)$.
\begin{theorem}\label{unitary_group_theorem} Let $c \in \Con(\R^d)$ satisfy
(\ref{coeff_structure}) and define $A v = \div( c(x) \grad v)$ with domain
$D(A) = H^2(\R^d)$ in $L^2(\R^d)$. Then $A$ is self-adjoint and $i\, A$
generates a strongly continuous unitary group $(T(z))_{z\in\R}$ on $L^2(\R^d)$.
Moreover, we have the following resolvent estimate, valid for $\la \in
\R\setminus\{0\}$:
\beq\label{resolvent_estimate}
   \norm{(iA - \la)^{-1}v}{L^2}
  \leq \frac{\norm{v}{L^2}}{|\la|}  \qquad \text{ for all } v \in L^2(\R^d).
\eeq
\end{theorem}

We briefly recall how $(T(z))_{z\in\R}$ can be used to construct solutions to
the Cauchy problem on $\R^d\times\R$
\[
  \d_z v - i A v = g, \quad v(0) = v_0.
\]
Let $v_0 \in L^2(\R^d)$ and $g \in L^1(\R,L^2(\R^d))$ then the \emph{mild
solution}
\beq\label{group_mild_solution}
    v(z) := T(z) v_0  + \int_0^z T(z - \rho) g(\rho) \, d\rho
\eeq
is in $\Con(\R,L^2(\R^d))$.

If $v_0 \in H^2(\R^d)$ and 
$g \in \Con(\R,H^2(\R^d))$ or
$g \in \Con^1(\R,L^2(\R^d))$  then $v$ belongs to  $\Con^1(\R,L^2(\R^d))$ 
and is the unique classical solution with pointwise values in $H^2(\R^d)$
(cf.\ \cite[Section 4.2, Corollaries 2.5 and 2.6]{Pazy:83}).

\begin{remark}\label{rem_mild_weak}
The mild solution (\ref{group_mild_solution})  defines a weak
solution in the following sense: $v(0) = v_0$ and  for all
$\phi \in \D(\R^{d+1})$
\[
    - \int_\R \Big(
    \inp{v(z)}{\d_z \phi(z,.)- i A \phi(z,.)} \Big)\, dz
    = \int_\R \inp{g(z)}{\phi(z,.)}\, dz,
\]
where $\inp{\ }{\ }$ denotes the inner product in $L^2(\R^d)$. To see this, one
approximates the mild solution by classical solutions $(v_k)_{k\in\N}$ to
equations with regularized right-hand side and initial data (\cite[Section 4.2,
Theorem 2.7]{Pazy:83}):  $L^2$-convergence of $v_k(z) \to v(z)$ (as $k \to
\infty$), uniformly when $z$ varies in compact intervals, together with the
convergence $g_k \to g$ in $L^1(\R,L^2(\R^d))$ implies convergence in the
integral formula above.
\end{remark}

\subsection{Evolution system at fixed frequency $\tau$}

Let $\tau\in\R$ be fixed, but arbitrary. We consider the $z$-parameterized
family of unbounded self-adjoint operators in (\ref{A_operator}) and put
 \[
    A(\tau;z) := A(\tau;z,x,D_x) \quad (z \geq 0).
 \]
Let $Z > 0$ be arbitrary. We will check that, for every $\tau$, $(i
A(\tau;z))_{z \geq 0}$ defines an evolution system (or fundamental solution)
$(U(\tau;z_1,z_2))_{Z \geq z_1 \geq z_2 \geq 0}$ on $L^2(\R^d)$ by applying
\cite[Section 4.4, Corollary to Theorem 4.4.2, p.\ 102]{Tanabe:79} (cf.\ also
\cite[Sections 5.3-5.5]{Pazy:83}). We have to check that the corresponding
hypotheses are satisfied.

First, observe that $D(A(\tau;z)) = H^2(\R^d)$ is independent of the evolution
parameter $z$ (and of $\tau$), and every $i A(\tau;z)$ is the skew-adjoint
generator of a strongly continuous (unitary) semigroup $(T(\tau,z;\zeta)_{\zeta
\geq 0})$ on $L^2(\R^d)$. Furthermore, the resolvent estimates
(\ref{resolvent_estimate}), valid for all $z$ (and $\tau$), immediately imply
that $(A(\tau;z))_{z \geq 0}$ is a stable family of generators \emph{with
stability constants $1$ and $0$ for all $\tau$} (cf.\ \cite[Definition
4.3.1]{Tanabe:79}).

Finally, we have to check that for all $v \in H^2(\R^d)$ the map
 \[
    [0,\infty) \ni z \mapsto A(\tau;z) v \in L^2(\R^d)
 \]
 is continuously differentiable. We may use equation (C2) from Fact C to write
 (with $\grad$ taken with respect to $x$ only)
  \[
    A(\tau;z) v = \grad c(z,x,\tau) \cdot \grad v +  c(z,x,\tau) \, \Delta v.
  \]
  By Assumption 1,(ii-iii), we have
  \[
        \grad c(.,.,\tau) \in \Con^1([0,\infty),H^r(\R^d)), \quad
            c(.,.,\tau) - c_0 h_0(\tau) \in \Con^1([0,\infty), H^{r+1}(\R^d)).
  \]
Since $\grad v \in H^1(\R^d)$ and $\Delta v \in L^2(\R^d)$ the multiplication
rules plus continuity properties in Fact A apply (where in case $H^r \cdot H^1$
we choose $\eps < r$, if $d=2$) and yield that $A(\tau;.) v \in
\Con^1([0,\infty),L^2(\R^d))$.

Thus, all hypotheses of \cite[Section 4.4, Corollary to Theorem 4.4.2, p.\
102]{Tanabe:79} are fulfilled. Note that the evolution system is constructed as
the strong operator limit of discretizations based on the unitary semigroups of
each generator, hence is contractive. This implies the
following intermediate result.
\begin{proposition}\label{evolution_system_proposition} Let $Z > 0$. Then
for all $\tau\in\R$ the family $(i A(\tau;z))_{z \geq 0}$ defines a unique
evolution system $(U(\tau;z_1,z_2))_{Z \geq z_1 \geq z_2 \geq 0}$ on
$L^2(\R^2)$ with the following properties: The map $(z_1,z_2) \mapsto
U(\tau;z_1,z_2)$ is strongly continuous, $U(\tau;z,z) = I$, $U(\tau;z_1,z_2)$
is contractive, and
 \beq\label{evolution_property}
    U(\tau;z_1,z_2) \circ U(\tau;z_2,z_3) = U(\tau;z_1,z_3) \qquad
        0 \leq z_3 \leq z_2 \leq z_1 \leq Z;
 \eeq
 moreover, $H^2(\R^d)$ is invariant under $U(\tau;z_1,z_2)$, for all $v \in
 H^2(\R^d)$  the map $(z_1,z_2) \mapsto U(\tau;z_1,z_2) v $ is 
 continuously differentiable, separately in both variables, and the 
following equations hold:
 \begin{align}
    \pdiff{z_1} U(\tau;z_1,z_2) v & = A(\tau;z_1) U(\tau;z_1,z_2) v \\
    \pdiff{z_2} U(\tau;z_1,z_2) v & = - U(\tau;z_1,z_2) A(\tau;z_2) v.
 \end{align}
\end{proposition}

At this stage, we obtain solutions to a version of the Cauchy problem
(\ref{PDE}-\ref{initial_cond}) at fixed frequency $\tau$, i.e.,
\begin{align}
  \d_z v - i A(\tau;z) v & = g \in L^1([0,Z],L^2(\R^d)) \label{PDE_tau}\\
    v \mid_{z=0}  & = v_0 \in L^2(\R^d). \label{initial_cond_tau}
\end{align}
%If $v_0 \in L^2(\R^d)$ and $g \in L^1([0,\infty),L^2(\R^d))$ then
The \emph{mild solution} is defined by
\beq\label{mild_solution_tau}
    v(z) := U(\tau;z,0) v_0  + \int_0^z U(\tau;z,\rho) g(\rho) \, d\rho
\eeq and belongs to $\Con([0,Z],L^2(\R^d))$ (\cite[Section 5.5, Definition
5.1]{Pazy:83}).

\begin{remark}\label{rem_regularity}
(i) In the case of classical solutions, we have the following regularity property:
If $v_0 \in H^2(\R^d)$ and $g \in \Con([0,Z],H^2(\R^d))$ or $g \in
\Con^1([0,Z],L^2(\R^d))$  then $v \in \Con^1([0,Z],L^2(\R^d))$ 
is the unique \emph{$H^2$-valued solution} and 
satisfies the equation in the strong sense (cf.\ \cite[Section 5.5, Theorems
5.2 and 5.3]{Pazy:83}).

(ii) Observe that, at frozen value of $\tau$, one may 
  apply \cite[Chapter 3, Theorem 10.1 and Remark 10.2]{LM:72} directly
  by putting $H = L^2$, $V = H^1$, and
\[
        a(z;u,v) = \sum_{j=1}^{d} \inp{c(z,.) \d_j u}{\d_j v}
                \qquad \forall u, v \in V.
\]
It suffices to assume $c \in \Con^1([0,Z],L^\infty(\R^d))$, then for
any initial value $v_0 \in H^1(\R^d)$ and right-hand side $g \in
L^2([0,Z]\times\R^d)$ such that $\d_z g \in L^2([0,Z],H^{-1})$ there
is a unique solution $v \in \Con([0,Z],H^1) \cap \Con^1([0,Z],H^{-1})$
to the Cauchy problem (\ref{PDE_tau}-\ref{initial_cond_tau}). However,
our approach allows for a precise investigation of the
$\tau$-dependence, which is needed to solve the full Cauchy problem
(\ref{PDE}-\ref{initial_cond}) with distributional data as well as to
transform back to the original problem
(\ref{orig_PDE}-\ref{orig_initial_cond}) in Section 4. Furthermore,
our results show that lateral $H^2$-regularity of the data is
preserved in the solution.
 \end{remark}

We thus have established an evolution system in the $L^2$-setting.
Note that by Lemma \ref{elliptic_regularity_prop} we have, in fact,
that $A(\tau;z,x,D_x)$ is an unbounded operator on $H^s$ with domain
$H^{s+2}$ for any $0 \leq s < r$.  If we were able to establish an
evolution system on $H^s$ then the regularity information encoded into
$A$ would be more directly preserved. 

\subsection{Frequency dependence of the evolution system}

Throughout this subsection, let $Z > 0$ be arbitrary but fixed. So far, the
frequency parameter $\tau$ was arbitrary, but fixed, throughout the
construction of the evolution system $(U(\tau;z_1,z_2))_{Z \geq z_1 \geq z_2
\geq 0}$. We will prove that the dependence on all parameters $(\tau,z_1,z_2)$
jointly is strongly continuous. In the sequel, let $L(E,F)$ (resp.\ $L(E)$)
denote the set of bounded linear operators between the Banach spaces $E$ and
$F$ (resp.\ on $E$).

We begin with an observation on the general level of semigroups and evolution
systems.
\begin{lemma}\label{U_continuity} Assume that
\beq\label{A_continuity_condition}
  (\tau,z) \mapsto A(\tau;z) \text{ is continuous }
     \R\times [0,\infty) \to L(H^2(\R^d),L^2(\R^d))
\eeq
(with respect to the operator norm) and
\beq\label{T_continuity_condition}
  (\tau,z,\zeta) \mapsto T(\tau,z;\zeta) \text{ is strongly continuous }
     \R\times [0,\infty)\times [0,\infty) \to L(L^2(\R^d)),
\eeq where $(T(\tau,z;\zeta))_{\zeta \geq 0}$ denotes the semi-group generated
by $A(\tau;z)$. Then the map $(\tau,z_1,z_2) \mapsto U(\tau;z_1,z_2)$ is
strongly continuous from $ B:= \R \times \{ (z_1,z_2) \colon Z \geq z_1 \geq
z_2 \geq 0 \}$ into $L(L^2(\R^d))$.
\end{lemma}
\begin{proof}
We inspect the basic construction of the evolution system from the family of
semigroups in the proof of \cite[Section 5.3, Theorem 3.1]{Pazy:83} and keep
track of the additional parameter $\tau\in\R$ in our case. For all
$(\tau,z_1,z_2) \in B$  we obtain $U(\tau;z_1,z_2)$ as the strong limit
 of $U_n(\tau;z_1,z_2)$ (as $n\to\infty$), where $U_n(\tau;.,.)$ is the evolution
 system defined as follows: put $z_n^j = j Z / n$ ($j=0,\ldots,n$)
 then for $\tau\in\R$, $0 \leq y \leq z \leq Z$ let
 $$ 
    U_n(\tau;z,y) := T(\tau,z_n^l;z-y)
        \qquad\qquad
            \text{if } z_n^l \leq y \leq z \leq z_n^{l+1},
$$ 
and
\begin{multline*}
    U_n(\tau;z,y)  := T(\tau,z_n^k;z-z_n^k) \cdot\!\!\!\!\!\!
                        \prod\limits_{l+1 \leq j \leq k-1} \!\!\!\!\!\!
                             T(\tau,z_n^j;Z/n)
                        \cdot T(\tau,z_n^l;z_n^{l+1}-y)\\
                    \text{if } z_n^l \leq y \leq z_n^{l+1} \leq z_n^k
                            \leq z \leq z_n^{k+1}, k > l.
\end{multline*}
By (\ref{T_continuity_condition}) the right-hand side of each formula is
strongly continuous with respect to $(\tau,z,y)$, and the boundary values,
when $k = l+1$ and $y = z_n^{l+1}$ or $z = z_n^{l+1}$, match. Hence $U_n$ is
strongly continuous on $B$ and $\| U_n(\tau;z,y) \| = 1$.

As in \cite[(3.13) on p.\ 136]{Pazy:83} we have the following integral
representation for the action on any $v \in H^2$
 \[
    U_n(\tau;z,y) v - U_m(\tau;z,y) v = \int_y^z U_n(\tau;z,\rho)
        \big( A_n(\tau;\rho) - A_m(\tau;\rho) \big)
            U_m(\tau;\rho,y) v \, d\rho,
 \]
where $A_n(\tau;\rho)$ is the piecewise constant approximation of
$A(\tau;\rho)$ with $A_n(\tau;\rho) := A(\tau;z_n^k)$, when $z_n^k \leq \rho <
z_n^{k+1}$, and $A_n(\tau;Z) = A(\tau;Z)$. By (\ref{A_continuity_condition}) we
have $\norm{A_n(\tau;\rho) - A(\tau;\rho)}{L(H^2,L^2)} \to 0$ uniformly for
$(\tau,\rho)$ in compact sets. Passing to the limit $m \to \infty$ in
the integral 
representation above yields the estimate
\[
    \ltwo{U_n(\tau;z,y) v - U(\tau;z,y) v} \leq \norm{v}{H^2} \,
        \int_y^z  \norm{A_n(\tau;\rho) - A(\tau;\rho)}{L(H^2,L^2)}\, d\rho.
\]
By the uniform convergence of $A_n(\tau;\rho)$ (as $n \to \infty$) we thus
obtain (local) uniform convergence of $U_n(\tau;z,y)v$, which proves the
asserted continuity of $(\tau,z,y) \mapsto U(\tau;z,y) v$.
\end{proof}

We have to establish conditions
(\ref{A_continuity_condition}-\ref{T_continuity_condition}) in the specific
context of the assumptions described in Section 2. In due course, we will make
repeated use of $(\tau,z)$-parameterized variants of Facts A-C, stated in
Subsection 3.1. Note that, in particular, the function $F$ used in Fact B does
not depend on $(\tau,z)$.
\begin{lemma}\label{A_continuity} If $A(\tau;z)$ ($\tau\in\R$, $z \geq 0$) is
defined by (\ref{A_operator}) then Assumption 1 implies Lipschitz-continuity of
the map in condition (\ref{A_continuity_condition}).
\end{lemma}
\begin{proof}
Let $M(\tau,z)$ denote the operator of multiplication of pairs $(v_1,v_2) \in
H^1 \times H^1$ by the scalar function $c(\tau,z,x) - c(\tau_0,z_0,x)$. We
write $A(\tau;z) - A(\tau_0;z_0) = \div \circ M(\tau,z) \circ \grad $ as a
composition of operators and get the following norm inequality
  \begin{multline*}
     \norm{A(\tau;z) - A(\tau_0;z_0)}{L(H^2,L^2)} \\ \leq
        \norm{\div}{L(H^1 \times H^1,L^2)} \cdot \norm{M(\tau,z)}{L(H^1 \times H^1)}
         \cdot \norm{\grad}{L(H^2,H^1 \times H^1)} \\
         \leq \sqrt{2}\, \norm{M(\tau,z)}{L(H^1 \times H^1)}.
  \end{multline*}
To estimate $\norm{M(\tau,z) (v_1,v_2)}{H^1 \times H^1}$ it suffices to find an
upper bound of \\ $\norm{(c(z,.,\tau) - c(z_0,.,\tau_0))\, v}{H^1}$ for $v\in
H^1$. We have
 \begin{multline*}
    (c(z,x,\tau) - c(z_0,x,\tau_0)) \, v(x) = \\
        \begin{pmatrix} z - z_0 \cr \tau - \tau_0 \end{pmatrix}
        \cdot \int_0^1
        \grad_{(z,\tau)} c(z_0 + \sig (z - z_0), x, \tau_0 + \sig (\tau - \tau_0))
        \, d\sig \,\, v(x),
 \end{multline*}
 which, upon taking the $H^1$-norm with respect to $x$ and assuming
 $\max(|z - z_0|,|\tau - \tau_0|) \leq 1$, yields
  \begin{multline*}
    \norm{(c(z,.,\tau) - c(z_0,.,\tau_0))\, v}{H^1} \leq
        \max (|z - z_0|, |\tau - \tau_0|)  \\ \cdot
        \sup (\norm{\d_z c(z',.,\tau') v}{H^1} + \norm{\d_\tau
        c(z',.,\tau')}{H^1}),
  \end{multline*}
 where the supremum is taken over
 $(z',\tau') \in [z_0-1, z_0+1]\times[\tau_0-1,\tau_0+1]$.
Assumption 1 implies that $\d_z c$, $\d_\tau c$ both are continuous functions
of $(z,\tau)$ valued in $H^{r+1}(\R^d)$, which combined with Fact A gives
 \begin{multline*}
    \norm{\d_z c(z',.,\tau') v}{H^1} + \norm{\d_\tau c(z',.,\tau')}{H^1} \\
    \leq C_1 \, \norm{v}{H^1} \,
        (\norm{\d_z c(z',.,\tau')}{H^{r+1}} + \norm{\d_\tau
        c(z',.,\tau')}{H^{r+1}})\leq
        C_2 \, \norm{v}{H^1}
 \end{multline*}
 with positive constants $C_1$, $C_2$ and for all
 $(z',\tau') \in [z_0-1, z_0+1]\times[\tau_0-1,\tau_0+1]$. Combining all estimates
 we deduce that there is $C_3
  > 0$ such that $|z - z_0| + |\tau - \tau_0| \leq 1$ implies
  \[
     \norm{M(\tau,z) (v_1,v_2)}{H^1 \times H^1} \leq C_3 \,
        \max (|z - z_0|, |\tau - \tau_0|)\,  \norm{(v_1,v_2)}{H^1 \times H^1},
  \]
which proves the asserted Lipschitz-continuity.
\end{proof}

\begin{lemma}\label{T_continuity} If $A(\tau;z)$ ($\tau\in\R$, $z \geq 0$) is
defined by (\ref{A_operator}) then Assumption 1 implies the continuity
condition (\ref{T_continuity_condition}).
\end{lemma}
\begin{proof}
We apply the Kato-Trotter theorem on
convergence of semi-groups 
(cf.\ \cite[Chapter IX, Section 12, Theorem 1]{Yosida:80}).
According to this theorem, we obtain
\begin{center}  $T(\tau,z;\zeta) \to T(\tau_0,z_0;\zeta)$ strongly as
        $(\tau,z) \to (\tau_0,z_0)$, \\
    uniformly on any compact interval containing $\zeta$,
\end{center}
thus (\ref{T_continuity_condition}) by uniformity, provided that we show strong
continuity of the resolvent map $(\tau,z) \mapsto (\la - i A(\tau;z))\inv =:
R(\la,iA(\tau;z))$ for some $\la > 0$.

Fix $\la > 0$ and let $f\in L^2(\R^d)$ be arbitrary. Define $u(\tau,z) :=
R(\la,iA(\tau;z))f \in H^2(\R^d)$, so that $u$ solves
 \beq\label{resolvent_equation}
    \la \, u(\tau,z) - i A(\tau;z) u(\tau,z) = f.
 \eeq
Adding the difference $i A(\tau;z)u(\tau,z) - i A(\tau_0;z_0)u(\tau,z)$ yields
\[
    (\la - i A(\tau_0;z_0))\, u(\tau,z) =
        f + i (A(\tau;z)- A(\tau_0;z_0))\, u(\tau,z) =: f + i w(\tau,z).
\]
Hence, $u(\tau,z) = R(\la,i A(\tau_0;z_0)) (f + i w(\tau,z))$ and it suffices
to prove that $w(\tau,z) \to 0$ in $L^2(\R^d)$ as $(\tau,z) \to (\tau_0,z_0)$.
Applying (C2) from Fact C in Subsection 3.1, we may write
\begin{multline*}
    w(\tau,z) = \grad \big(c(z,x,\tau)-c(z_0,x,\tau_0) \big)
        \cdot \grad u(\tau,z) \\
    + \big( c(z,x,\tau)-c(z_0,x,\tau_0) \big)\, \Delta u(\tau,z).
\end{multline*}
By Assumption 1, the difference $c(z,.,\tau)-c(z_0,.,\tau_0)$ tends to $0$ in
$H^{r+1}(\R^d)$. In view of Fact A this implies $w(\tau,z) \to 0$, if $\grad
u(\tau,z)$ as well as $\Delta u(\tau,z)$ stays bounded. To prove the latter, we
take the $L^2$-inner product with $u$ on both sides of equation
(\ref{resolvent_equation}) and obtain
 \[
    \la \, \ltwo{u}^2 - i \sum_j \inp{c\, \d_{x_j}u}{\d_{x_j}u} = \inp{f}{u}.
 \]
 Note that taking real parts here yields the estimate
 (\ref{resolvent_estimate}), which is $\ltwo{u} \leq \ltwo{f}/\la$ . If we
 take absolute values of the imaginary parts, we may use the lower bound
 $c(z,x,\tau) \geq c_0$ and the resolvent estimate to deduce $\sum_j
 \ltwo{\d_{x_j}u(\tau,z)}^2 \leq \ltwo{f}^2 / (\la c_0)$, uniformly in
 $(\tau,z)$. 

Finally, the boundedness of $\ltwo{\Delta u(\tau,z)}$ is revealed in several
steps. 
First, note that (C2) from Fact C applied to (\ref{resolvent_equation})
yields 
\[
        c(z,.,\tau)\, \Delta u(\tau,z) = i(f - \la\, u(\tau,z)) - 
        \grad c(z,.,\tau) \cdot \grad u(\tau,z).
\]
The first term on the right-hand side is bounded in $L^2$, 
uniformly for all $(\tau,z)$, whereas the second term is uniformly 
bounded in $H^{r-1}$ by Fact A. Hence $c\,\Delta u(\tau,z)$ is a bounded 
family in $H^{r-1}$ and, combining Facts A and C, we find that $\Delta
u(\tau,z)$ is uniformly bounded in $H^{r-1}$ as well. Therefore, $u(\tau,z)$
has a uniform bound in $H^{r+1}$-norms. From here we may proceed as in Claim 3
from the proof of Proposition \ref{elliptic_regularity_prop} (with $s = 0$).
Indeed, the arguments used there preserve uniform boundedness properties
throughout, since we have such in $H^{r+1}$ already. Thus, $u(\tau,z)$ is
uniformly bounded in $H^2$, in particular, $\Delta u(\tau,z)$ is bounded
uniformly for all $(\tau,z)$, which completes the proof.
\end{proof}

We summarize the preceding results in the announced continuity statement for
the evolution system.
\begin{theorem} Let $(U(\tau;z_1,z_2))_{Z \geq z_1 \geq z_2 \geq 0}$ be the
evolution system generated by the family of operators $A(\tau;z)$
($\tau\in\R$, $z \geq 0$), defined in (\ref{A_operator}) and satisfying
Assumption 1. Then $(\tau,z_1,z_2) \mapsto U(\tau;z_1,z_2)$ is strongly
continuous $\R \times \{ (z_1,z_2) \colon Z \geq z_1 \geq z_2 \geq 0 \} \to
L(L^2(\R^d))$.
\end{theorem}

%%% Local Variables: 
%%% mode: latex
%%% TeX-master: "paraxial"
%%% End: 

\section{Solution of the Cauchy problem}

In this section we present our main results: existence and
uniqueness of solutions to the Cauchy problem 
(\ref{PDE}-\ref{initial_cond}) in the frequency domain and to
(\ref{orig_PDE}-\ref{orig_initial_cond}) in the time domain.

If $E$ is a Banach space, let $\Con_b(\R,E)$ denote the space of
$E$-valued continuous bounded functions.
Observe that for any $G\in \Con_b(\R,L^2) \subset \S'(\R;L^2)$ the expression
$U(\tau;z_1,z_2) G(\tau)$ is well-defined pointwise for all
$(\tau,z_1,z_2)$. Therefore, collecting the results obtained so far we
arrive at the following assertion.
\begin{proposition} If $v_0\in \Con_b(\R,L^2(\R^d))$ and
$g\in\Con_b([0,Z]\times\R,L^2(\R^d))$ 
then the formula
\beq\label{mild_solution}
    v(z,\tau) := U(\tau;z,0) v_0(\tau)  +
        \int_0^z U(\tau;z,\rho) g(\rho,\tau) \, d\rho
\eeq
defines a mild solution 
$v\in\Con^1([0,Z],\Con_b(\R,L^2(\R^d)))\subset \SW$ to
(\ref{PDE}-\ref{initial_cond}). Moreover, when $v$ is a strong
solution then $u := \F_t\inv v$ is a strong solution of 
(\ref{orig_PDE}-\ref{orig_initial_cond})
with initial data $u_0 = \F_t\inv v_0$ and right-hand side 
$f = \F_t\inv g$. 
\end{proposition}
For example, the hypotheses leading to strong solvability
 are satisfied if
$v_0 \in \Con_b(\R,H^2)$ and $g \in \Con_b([0,Z]\times\R,H^2)$ or
$g \in \Con^1([0,Z],\Con_b(\R,L^2))$. Of course, using
functions that are bounded and continuous with respect to the frequency variable
$\tau$ here is just one simple way to ensure that all constructions
described above work with all involved objects staying temperate. More
generally, it would suffice to consider elements in $\SW$ whose
distributional action (with respect to the frequency variable) is
given by (weak) integration over a continuous function (times the test
function). 

To apply formula (\ref{mild_solution}) to the original Cauchy problem
(\ref{orig_PDE}-\ref{orig_initial_cond}) we only need to state
conditions on the data $u_0$ and $f$ that imply $\F_t u_0 \in
\Con_b(\R,H^2)$ and $\F_t f \in \Con_b([0,Z]\times\R,H^2)$ or $\F_t f
\in \Con^1([0,Z],\Con_b(\R,L^2))$. Note that, for example, in our
physical application
such conditions would be met if the source or force terms are active
only for some finite time interval and vanish otherwise.
Then by the uniqueness of
$H^2$-valued solutions, as stated in Remark \ref{rem_regularity}, we obtain the
following result.

\begin{theorem}\label{summary_thm}
Assume that the right-hand side $f$ in equation (\ref{orig_PDE}) satisfies either 
$f \in \Con\big([0,Z],L^1(\R,H^2(\R^2))\big)$
or $f \in \Con^1\big([0,Z],L^1(\R,L^2(\R^2))\big)$.

Then for every $u_0 \in L^1(\R,H^2(\R^2))$  
the Cauchy problem
(\ref{orig_PDE}-\ref{orig_initial_cond}) has a unique
strong solution 
$u\in\Con^1([0,Z],\S'(\R,L^2(\R^2)))$ which is $H^2$-valued in the 
following sense: for all $z\in [0,Z]$
and $\phi \in \S'(\R)$ we have  $\dis{u(z)}{\phi} \in H^2(\R^2)$.

Moreover,  $u$ is obtained by inverse partial Fourier transform (with
respect to $\tau$) of
$v\in\Con^1([0,Z],\Con_b(\R,L^2(\R^2)))$ as defined in Equation
(\ref{mild_solution}), where $v_0 := \F_t u_0$, $g := \F_t f$, and $U$
is the evolution system from Proposition
\ref{evolution_system_proposition}. In addition,  
$v(z,\tau)$ belongs to $H^2(\R^2)$ for every $(z,\tau)$.
\end{theorem}

\paragraph{Inverse analysis of medium regularity:}
We conclude with a brief indication of a potential application of 
Theorem \ref{summary_thm} to an inverse analysis of medium
regularity in wave propagation. Suppose that we are in a model situation
where parts (i-ii) and (iv) of Assumption 1 in Section 2 are satisfied and the
regularity property (iii) of the medium is in question; assume that we
have the a priori knowledge that $c_l \in
\Con^1([0,\infty),L^\infty(\R^2))$ 
for all $l$. Let the sources of a seismic experiment be calibrated to 
produce data in accordance with the hypotheses in Theorem
\ref{summary_thm}. If then the measured wave solution $u$ (or $v$) fails
to display the asserted $H^2$-regularity then we may conclude
that for some $c_l$ (and near some depth $z$) there is no $r \in (0,1)$
such that the $H^{r+1}$-regularity holds; in other words, the
lateral regularity of the medium there cannot be better than
$H^1$ (on the Sobolev scale). If we were in the possession of 
analogous $H^s$-results ($s > r$)  
for the Cauchy problem it would enable us to draw sharper 
conclusions in such a ``inverse regularity analysis''. Note that the
exact location of the most singular region need not be known. For the
application hence, precise imaging of the singularities is not
required prior to the regularity analysis.

%%% Local Variables: 
%%% mode: latex
%%% TeX-master: "paraxial"
%%% TeX-master: "paraxial"
%%% End: 

\bibliographystyle{abbrv}
\bibliography{gue}

\end{document}